\documentclass[12pt]{amsart}
\usepackage{amscd}
\usepackage{verbatim}
\usepackage{amssymb, amsmath, amsthm, amscd,ifthen}
\usepackage[dvips]{graphics}
\usepackage[cp866]{inputenc}
\usepackage{graphicx}
\usepackage{epsfig}
\usepackage{color}
\usepackage{tikz-cd}
\usepackage{amsmath,amssymb,amscd}
\usepackage[all]{xy}
\usepackage{mathtools}

\usepackage{amssymb}

\newtheorem{thm}{Theorem}[section]

\newtheorem{prop}[thm]{Proposition}

\newtheorem{cor}[thm]{Corollary}

\theoremstyle{definition}

\newtheorem{dfn}[thm]{Definition}

\newcommand{\M} {\mathcal M}

\newcommand{ \HH }{H}
\newcommand{ \ZZ } {\mathbb{Z}}

\DeclareMathOperator{\UU}{U}

\title[On the cohomology ring of polygon spaces]{Geometric presentation for the cohomology ring of polygon spaces}
\author{Ilia Nekrasov}

\address{St.~Petersburg State University, Chebyshev Laboratory}
  \email{geometr.nekrasov@yandex.ru}

\author{Gaiane Panina}

\address{St.~Petersburg Department of V.A. Steklov Institute of Mathematics of the Russian Academy of Sciences}
\email{gaiane-panina@rambler.ru}


\begin{document}

\begin{abstract}
We describe the cohomology ring of the moduli space of a flexible polygon in geometrically meaningful terms. We  propose two presentations, both are  computation friendly: there are simple rules for cup product.

\end{abstract}

\maketitle \setcounter{section}{0}

\section{Introduction}

With an $n$-tuple of positive numbers $L=(l_1,...,l_n)$ we associate a \textit{flexible polygon}, that is, $n$  rigid bars of lengths $l_i$ connected in a   cyclic chain by revolving joints. A \textit{configuration} of   $L$ is an $n$-tuple  of points $(q_1,...,q_n)\in (\mathbb{R}^3)^n$ with $|q_iq_{i+1}|=l_i$ and $ |q_nq_1|=l_n$.

\begin{dfn}\label{defConfSpace}
 \textit{The moduli space
 $\M_3(L)$}  is the set  of all configurations  of $L$ lying in $\mathbb{R}^3$
modulo orientation preserving  isometries of $\mathbb{R}^3$.
\end{dfn}
\begin{dfn}\label{SecondDfn}
Equivalently, one defines

$$\M_3(L)=\left\{(u_1,...,u_n) \in (\mathbb{S}^2)^n : \sum_{i=1}^n l_iu_i=0\right\}/SO(3). $$
\end{dfn}

The second definition shows that   $\M_{3}(L)$ does not depend on the
ordering of $\{l_1,...,l_n\}$; however, it depends on the values
of $l_i$.

The space $\M_{3}(L)$  is a $2n-6$-dimensional complex-analytic manifold.\footnote{Moreover, A.~Klyachko  \cite{Kl} showed that it is an algebraic variety.} The cohomology ring $\HH^{*}(\M_{3}(L), \ZZ)$   is explicitly described in \cite{HKn}. Some of its generators and relations have a clear geometrical meaning, while the other do not.
In \cite{Kami} Kamiyama writes: "Although (R1) and (R2) are manageable,  (R3) is complicated and mysterious. It is not easy to obtain information from (R3)."\footnote{{In notation of the present paper, (R1), (R2), and (R3) are relations (2), (3), and (4) from Theorem \ref{CR_HKn} respectively.}} In the cited paper he gives a reasonable answer by providing a \textit{combinatorial comprehension} of the relation modulo $\mathbb{Z}_2$.

In turn, the present paper interprets generators and relations in entirely geometric terms. The \textbf{first presentation } (that is, the ring $\widetilde{\mathcal{N}}$) is described in Section \ref{computation}.
Oversimplifying, its elements are $\mathbb{Z}$--linear combinations of (Poincar\'{e} duals of)\textit{ nice submanifolds} of $\M_{3}(L)$. The latter are characterized  by codirected and oppositely directed pairs of edges, see Section \ref{nicemanifolds}.

The first and the second series of relations in the ring $\widetilde{\mathcal{N}}$ have a transparent geometrical meaning.

The third series of relations comes from computation of Chern classes of some natural bundles over the moduli space.
The computation has in the background the (scholarly exercise!) way of finding the Chern class of tangent bundle $T\mathbb{S}^2$ to the $2$-sphere $\mathbb{S}^2$, see comments to Proposition \ref{Ch}, and also \cite{NPZ} for more details.

These bundles are quite remarkable: on the one hand,  they naturally appear in \cite{HKn}. On the other hand, they are a counterpart of M. Kontsevich's tautological bundles over the  space $\M_{0,n}$  \cite{Kon}.

In the  \textbf{second presentation } (the ring $\mathcal{P}$, see Section \ref{SecPerf}) the generators are (Poincar\'{e} duals of)\textit{ perfect submanifolds} of $\M_{3}(L)$. The latter are characterized by codirected edges only.


\textbf{Acknowledgement.} This research is supported by the Russian Science Foundation under grant 16-11-10039.

\section{Preliminaries}\label{prelim}
\subsection*{Moduli space $\mathcal{M}_{3}(L)$}\label{modulispace}

Assume that an $n$-tuple of positive real numbers \newline $L = (l_{1}, \dots, l_{n})$ is fixed. The moduli space  $\M_{3}(L)$ is not empty if and only if $L$ satisfies the  inequalities:
$$l_{i} < \frac{l_{1}+ \dots +l_{n}}{2} \;\text{ for all }i \in[n] = \{1, 2, \dots, n\}.$$
A subset $I \subset [n]$ is called \textit{long} if$$\sum_Il_i>\frac{1}{2}\sum_{i=1}^{n}l_i.$$ Otherwise, $I$ is  called \textit{short}.

Throughout the paper we assume that no configuration  fits in a straight line. This happens iff the moduli space $\M_{3}(L)$ is a smooth closed manifold. In more details, let us take all subsets $I\subset \{1,...,n\}$. The associated hyperplanes
$$\sum_{i\in I}l_i =\sum_{i\notin I}l_i$$
called \textit{walls }subdivide the parameter space $\mathbb{R}_+^n$ into a number of \textit{chambers}. The diffeomorphic type of $\M_3(L)$  depends only on the chamber containing $L$; this becomes  clear in view of the (coming below) stable configurations representations. For $L$ lying strictly inside a chamber, the space $\M_{3}(L)$ is a smooth manifold.

Let us make an additional assumption: the sum $\sum_I\pm l_i$ never vanishes for all non-empty $I\subset \{1,...,n\}$. This agreement does not restrict generality: one may perturb the edge lengths while staying in the same chamber. 

\medskip

We make use of yet another representation of   $\M_{3}(L)$.
Consider  configurations of $n$ (not necessarily all distinct)
points
$p_i$ in the complex projective line.  Each point $p_i$  is assigned the weight
 $l_i$.  The
configuration of (weighted) points is called  {\em
stable} if sum of the weights of coinciding points is
less than half the weight of all points.
Denote by  $S_\mathbb{C}(L)$ the space of stable configurations in the complex projective line.
The group $PSL(2,\mathbb{C})$  acts naturally on this space.
In this setting we have:

$$\mathcal{M}_3(L)=S_\mathbb{C}(L)/PSL(2,\mathbb{C}).$$
  That is, for each $n$ we have a  finite series $\M_{3}(L)$ of compactifications of the space ${\mathcal{M}}_{0,n}(\mathbb{C})$ depending on the particular choice of the lengths $L$.

\medskip

\begin{thm}[\cite{HKn}]\label{CR_HKn} The (co)homology groups  of $\M_{3}(L)$ are free abelian.
The cohomology ring $\HH^{*}(\M_{3}(L), \ZZ)$ is isomorphic to the quotient of the ring $\ZZ[R, V_{1}, \ldots, V_{n-1}, U_{1}, \ldots, U_{n-1}]$ by the ideal $\mathcal{I}$ generated by the four families of elements:
\begin{enumerate}
\item $U_{i} - V_{i} - R$ for $i \in [n-1]$;

\item $U_{i}V_{i}$ for $i \in [n-1]$;

\item $\prod\limits_{i\in I}V_{i}$ for $I \subset [n-1]$ such that $I \cup \{n\}$ is long;

\item $$\sum_{S\subsetneq H} \left( \prod_{i \in S}V_{i} \right) \cdot R^{|H-S| -1}, $$
such that $H\subset [n-1]$ is a long set, and $S\subsetneq H$ ranges over all the sets such that  $S \cup\{n\}$ is short.

\end{enumerate}

\end{thm}

\subsection*{Nice manifolds \cite{NPZ}}\label{nicemanifolds}

Nice manifolds are  submanifolds of $\M_{3}(L)$. We define them as point sets, equip them  with orientation, and then focus on their  Poincar\'{e} dual cocycles.

\subsubsection*{Nice manifolds as point sets}

We start with an example:
let $i\neq j \in [n]$. Denote by $(ij)$  the image of the natural embedding of the space $\M_3(l_i+l_j,l_1,...,\hat{l}_i,...,\hat{l}_j,...,l_n)$ into
the space $\M_3(L)$.

That is, we think of the configurations of the new $n-1$-gon as the configurations of $L$ with parallel \textbf{codirected} edges $i$ and $j$ \textit{frozen} together to a single edge of the length $l_i+l_j$. Since the moduli space does not depend on the ordering of the edges, it is convenient to think that $i$ and $j$  are consecutive indices.

Analogously, we can define the space $(i\;\overline{j})$ of configurations with edges $i$ and $j$ parallel, but \textbf{oppositely directed}. In other words, the space $(i\;\overline{j})$ is an image of the space $\M_3(|l_i-l_j|,l_1,...,\hat{l}_i,...,\hat{l}_j,...,l_n)$ under natural embedding in the space $\M_3(L)$.

\medskip

Generally, we can freeze several collections of edges either oppositely directed or codirected, and analogously define a nice submanifold labeled by the formal product. All submanifolds arising this way are called \textit{nice submanifolds of $\M_3(L)$, } or just \textit{nice manifolds}  for short.

Putting the above more formally, each nice manifold is labeled by an unordered formal product $$(I_1\overline{J}_1)\cdot ...\cdot(I_k\overline{J}_k),$$ where $I_1,...,I_k, {J}_1,,,{J}_k$ are some pairwise disjoint subsets of $[n]$ such that each set $I_i\cup {J}_i$ has at least one element.\footnote{For further computations is convenient to define also nice manifolds with $I_i\cup {J}_i$  consisting of one element. That is,  we set $(1)=\M_3(L),\ \  (1)\cdot (23)=(23)$ etc. }
{The manifold $(I_1\overline{J}_1)\cdot ...\cdot(I_k\overline{J}_k) $ is a subset of $\M_3(L)$ defined by the conditions:}
\begin{enumerate}
\item[(1)] $i,j \in I_k$ implies  $u_i=u_j$,

\item[(2)] $i,j \in J_k$ implies  $u_i=u_j$, and

\item[(3)] $i \in I_k, j\in J_k$ implies  $u_i=-u_j$.

\end{enumerate}

\medskip

\subsubsection*{Orientation on nice manifolds}
By construction, each nice manifold is the moduli space $\M_{3}(\tilde{L})$
for some length vector $\tilde{L}$. Therefore a nice manifold is a complex-analytic manifold and has a \textit{canonical} orientation  coming from the complex structure.

{By definition, the  \textit{relative orientation} of a nice manifold $(I\overline{J})$  coincides with its canonical orientation iff $$\sum _I l_i>\sum _J l_i.\eqno{(\ddag)}$$}

Further,  the two orientations (canonical and relative ones) of a nice manifold$$(I_1\overline{J}_1)\cdot ...\cdot (I_k\overline{J}_k)$$  coincide iff the above inequality $(\ddag)$  fails  for even number of $(I_i\overline{J_i})$.

From now on, by $(I_1\overline{J}_1)\cdot ...\cdot (I_k\overline{J}_k) \in H^*(\mathcal{M}_3(L),\mathbb{Z})$  we mean  (the Poincar\'{e} dual of) the  nice manifold taken with its relative orientation. Some of the nice manifolds might be empty and thus represent the zero cocycles.

\subsection*{{Cup product rules for nice manifolds \cite{NPZ}}}

The following rules are valid for nice submanifolds of $\M_3(L)$.
\begin{enumerate}
\item The cup product is a commutative operation.
\item $(I\overline{J})=-(J\overline{I}).$
  \item If  the factors have no common entries, the cup-product equals the formal product, e.g.:
  $$(12) \smile (34)\cdot (5\overline{6})=(12)\cdot (34)\cdot (5\overline{6}).$$

  \item If $I_1\cap I_2=\{i\},\ \ I_1\cap J_2= \emptyset, \ \ I_2\cap J_1= \emptyset,\  \hbox{and} \ \ J_1\cap J_2= \emptyset$,  then
$$(I_1\overline{J}_1)\smile (I_2\overline{J}_2)=
   (I_1\cup I_2 \ \overline{J_1\cup J_2}).
 $$
Examples: $$(123)\smile (345)=(12345),$$$$(123)\smile (34\overline{5})=(1234\overline{5}).$$
 \item If $J_1\cap J_2=\{i\},\ \ I_1\cap J_2= \emptyset, \ \ I_2\cap J_1= \emptyset,\ \hbox{and} \ \ I_1\cap I_2= \emptyset,$  then
$$(I_1\overline{J}_1)\smile (I_2\overline{J}_2)=-
   (I_1\cup I_2 \ \overline{J_1\cup J_2}).
 $$Example:
  $$(12\overline{3})\smile (45\overline{3})=  -(1245\overline{3}).$$

\end{enumerate}
Let us comment shortly on the rules: (1) comes from the complex structure, (2) is valid just by definition, and
 (3)--(5) reflect the transversal intersection of nice manifolds.

\subsection*{Natural bundles over $\M_{3}(L)$}\label{Natural}
Consider the following natural bundles over $\M_{3}(L)$:
\begin{dfn}\cite{HKn}
For $1 \leq i \leq n$ let $A_i=A_{i}(L)$ be the $\UU(1)$-bundle over $\M_{3}(L)$ whose fiber over a point $(u_{1}, \ldots, u_{n})$ consists of all rotations of the configuration around the vector $u_{i}$ (the direction of the rotation agrees with the right-hand rule).
\end{dfn}

\begin{dfn}\cite{NPZ}
The line bundle
 $E_i=E_{i}(L)$ is the complex line bundle over the space $\mathcal{M}_3(L)$  whose fiber over a point $(u_1,...,u_n)\in (\mathbb{C}P^1)^n$ is the complex tangent line to the complex projective line $\mathbb{C}P^1$
at the point $u_i$.
\end{dfn}
An easy observation is:
\begin{prop} \label{bundles}
The bundle $A_{i}$ is a principal $\UU(1)$-bundle corresponding to the complex line bundle $E_{i}$.\qed
\end{prop}

Denote by $Ch(i)$ the first Chern class of $A_i$ (or, equivalently, of $E_{i}$).

\begin{prop}\cite{HKn}\label{chernHKn}
In notation of Theorem \ref{CR_HKn}, one has
$$Ch(i) =
\begin{cases}
2V_{i} +R& \text{ for }i\neq n,\\
-R& \text{ for }i=n.\qed
\end{cases}
$$
\end{prop}

On the other side, we proved the analogous result in terms of nice manifolds:

\begin{prop}\cite{Agapito}, \cite{NPZ} \label{Ch}
For any $j \neq i \in [n]$ we have
$$Ch(i) = (ij) -(i \;\overline{j}).\qed$$
\end{prop}

Let us reveal the geometric meaning of this proposition.
We are going to describe a  section of $E_{i}$ which is transversal to the zero section. Its zero locus (modulo orientation) is the desired Chern  class  $Ch(i)$.
Fix any $j\neq i$, and consider the following (discontinuous) section of  $E_{i}$: at a point $(u_1,...,u_n)\in \mathcal{M}_3(L)$   we take the unit vector in the complex tangent line at $u_i$ pointing in the direction of the shortest arc connecting $u_i$ and $u_j$.  The section is well-defined except for the points of $\mathcal{M}_3(L)$ with $u_i=\pm u_j$, that is, except for nice manifolds $(ij)$ and $(i\,\overline{j})$.
One cooks a continuous section with zeroes  at $(ij)$ and $(i\,\overline{j})$.
  A detailed analysis (which we omit here;  the reader is referred to \cite{NPZ} for details) shows
that the signs are exactly as is stated in the proposition.

\medskip

\begin{prop} \label{Ch2}
For any $i \neq j \neq k\in [n]$ we have
$$Ch(i) = (ij) +(ik)-(jk).$$
\end{prop}
Proof. On the one hand, this relation can be derived from computation rules, see proof of Theorem \ref{Thmiso1}. On the other hand, it is good to understand the geometrical meaning:
Take a stable configuration $(x_1,...,x_n)\in \mathcal{M}_3(L)$. Take the circle passing through $x_i,x_j$, and $x_k$. It is oriented  by the order $ijk$. Take the  vector lying in the tangent complex line to $x_i$  which is tangent to the circle and points in the direction of $x_j$. It gives rise to a  section of $E_{i}$  which is defined correctly whenever these three points are distinct. Therefore, $Ch(i) = A(ij) +B(ik)+C(jk)$ for some integer $A,B,C$. Detailed analysis specifies their values.
\qed

\begin{cor}\label{CorChern2} (The four-term relation)
{  For any distinct $i,j,k,l$ we have}
  $$(ij)+(kl)=(jk)+(il).\qed$$
\end{cor}

\medskip

The above rules do not exhaust all existing products of nice manifolds. The following technical proposition shows a way to compute cup product in two particular cases that are not covered in the rules.

\begin{prop}\label{power}
For any $i,j,k$, we have
\begin{align*}
(ij) \smile (ij) &= (ij) \smile Ch(i) = (ijk)-(ij\overline{k}),\\
(i\overline{j}) \smile (i\overline{j}) &= (i\overline{j}) \smile Ch(i) = - (i\overline{j}k) +(i\overline{jk}).
\end{align*}
\end{prop}
\textit{Proof of  Proposition \ref{power}.
}
Using $(ij)\smile(i\overline{j}) = 0$  and the rules $(1)$ -- $(5)$,  we get
\begin{align*}
(ij)\smile (ij) &= (ij) \smile [(ij)-(i\overline{j})]
= (ij) \smile Ch(i) =\\
&= (ij)\smile [(ik)-(i\overline{k})] = (ijk)-(ij\overline{k}), \hbox{ and}
\end{align*}
\begin{align*}
(i\overline{j})\smile (i\overline{j}) &= (i\overline{j}) \smile [(i\overline{j})-(ij)]
= - (i\overline{j}) \smile Ch(i) =\\
&= - (i\overline{j})\smile [(ik)-(i\overline{k})] = - (i\overline{j}k) + (i\overline{jk}).\qed
\end{align*}

\section{The ring of nice manifolds. The first presentation }\label{computation}

Let us denote by $\mathcal{N}$ the $\ZZ$-linear closure of the set of all nice manifolds in the ring $\HH^{*}(\M_{3}(L), \ZZ)$.

\begin{thm}\label{theorem1}\begin{enumerate}
                             \item $\mathcal{N}$ coincides with the ring $\HH^{*}(\M_{3}(L), \ZZ)$.
In particular, $\mathcal{N}$ is closed under cup product.
                             \item  In notation of Theorem \ref{CR_HKn}, we have:
                           \newline  $R=-(ni)+(n\overline{i}),  $
                             $U_i=(n\overline{i}),$
                             $V_i=(ni).$
                           \end{enumerate}

\end{thm}
\textit{Proof.} {Prove first that the abelian group $\mathcal{N}$ is closed under cup product, so the set $\mathcal{N}$  is a {subring} of $\HH^{*}(\M_{3}(L), \ZZ)$.  It suffices to prove that each  monomial in elementary nice manifolds (i.e. manifolds $(i j)$ or $(i\overline{j})$ for some $i, j$) is a nice manifold.} So let us consider a monomial of the form

$$\pm \prod (i_p\overline{j}_p)^{\alpha_{p}} \smile \prod (i_q j_q)^{\beta_{q}}.$$

We may assume that $\Sigma \alpha_p+\Sigma\beta_q \leq n-3 = \dim \M_{3}(L)$, otherwise we get zero. Computation rules together with repeated applications  of Proposition \ref{power}  do the job.

Here are
two examples of such computation:
\begin{align*}
(123)\smile(124)&= (12)^2 \smile (13) \smile (14)=\\
&= [(125) - (12\overline{5})] \smile (13) \smile (14)=\\
&=(12345)-(123\overline{4}5).\\
(12)^{4}\smile (56)&= (12)^2 \smile [(123) - (12\overline{3})]\smile (56) =\\
&= (12) \smile [(1234) - (12\overline{3}4) - (123\overline{4}) + (12\overline{34})] \smile (56) =\\
&=(123456)-(12\overline{3}456)-(123\overline{4}56)-(1234\overline{5}6)+(12\overline{3}\overline{4}56)-(123\overline{4}\overline{5}6)\\
&
+(12\overline{3}4\overline{5}6)-(12\overline{345}6).
\end{align*}

On the one hand, the ring $\mathcal{N}$ is contained in the cohomology ring. On the other hand,  the cohomology ring is generated by nice manifolds due to Theorem \ref{CR_HKn}, Proposition \ref{chernHKn}, and Proposition \ref{Ch}.
\qed

\medskip
Now we are ready to describe the first computation friendly presentation of the cohomology ring.
Define the commutative graded ring $\widetilde{\mathcal{N}}$  with the ring operation "$*$" as follows:

(\textbf{Generators})   The generators of $\widetilde{\mathcal{N}}$ are the labels of all existing nice manifolds.
That is, the generators bijectively correspond to  all unordered formal products $$(I_1\overline{J}_1)\cdot ...\cdot(I_k\overline{J}_k),$$ where $I_1,...,I_k, {J}_1,,,{J}_k$ are some pairwise disjoint subsets of $[n]$ such that each set $I_i\cup {J}_i$ has at least one element.

(\textbf{Relations}) The relations in $\widetilde{\mathcal{N}}$ split in three sets:\begin{enumerate}
                                      \item Relations coming from "computation rules", see Section \ref{prelim}:
                                      \begin{enumerate}

\item $(I\overline{J})+(J\overline{I})=0.$
  \item If  the factors have no common entries, the product equals the formal product, e.g.:
  $$(12)* (34)=(12)\cdot (34).$$

  \item If $I_1\cap I_2=\{i\},\ \ I_1\cap J_2= \emptyset, \ \ I_2\cap J_1= \emptyset,\  \hbox{and} \ \ J_1\cap J_2= \emptyset$,  then
$$(I_1\overline{J}_1)* (I_2\overline{J}_2)=
   (I_1\cup I_2 \ \overline{J_1\cup J_2}).
 $$

 \item If $J_1\cap J_2=\{i\},\ \ I_1\cap J_2= \emptyset, \ \ I_2\cap J_1= \emptyset,\ \hbox{and} \ \ I_1\cap I_2= \emptyset,$  then
$$(I_1\overline{J}_1)* (I_2\overline{J}_2)=-
   (I_1\cup I_2 \ \overline{J_1\cup J_2}).
 $$

\medskip

\textbf{Remark.} Note that the relation $(d)$ follows from  $(a)$ and $(c)$.

\medskip

\item $(ij)*(i\overline{j})=0$. The geometrical meaning of this relation is transparent: the nice manifolds $(ij)$ and $(i\overline{j})$ are disjoint.

\end{enumerate}
                                      \item  A relation comes from length constraints:

                                        $(I\overline{J})=0$  whenever one of the sets $I$ or $J$ is long.

                                      \item  A relations coming from Chern class computation:

                                      for any $i,j,k \in [n]$ such that $i\neq j, i\neq k$, we have $$(ij)-(i\overline{j})=(ik)-(i\overline{k}).$$

                                    \end{enumerate}

\begin{thm}\label{Thmiso1}
  The rings $\widetilde{\mathcal{N}}$ and $H^*(\mathcal{M}_3(L),\mathbb{Z})$ are isomorphic.
\end{thm}
\textit{
Proof.} Since all the defining relations of $\widetilde{\mathcal{N}}$ are valid in $\mathcal{N}=H^*(\mathcal{M}_3(L),\mathbb{Z})$, we have an epimorphism
$$\widetilde{\mathcal{N}}\xrightarrow{\phi} H^*(\mathcal{M}_3(L),\mathbb{Z})=\mathcal{N},$$
where $\phi$ maps each formal expression to the dual of the associated nice manifolds.
To prove that $\phi$ is bijective, it suffices to show that for any set $H \subset [n-1]$, we have
$$\sum_{S \subsetneq H} (Sn) * (-Ch(n))^{*\big(|H-S| -1\big)} = \sum_{S \subsetneq H} (Sn) * (-Ch(n))^{*\big(|H-S| -1\big)} = (H),$$
where in the first sum, as in Theorem \ref{CR_HKn}, $S\subsetneq H$ ranges over all the sets such that  $S \cup\{n\}$ is short, and by $Ch(n)\in \widetilde{\mathcal{N}}$ we mean $(nj) -(n \overline{j})$.

The first equality is trivial, so we prove the second equality.

The case $|H|=2$  gives a base of induction.
 Without loss of generality, assume that $H = \{1,2\}$. The left-hand side  equals
\begin{align*}
-Ch(n) + (1n) +(2n)&= \frac{(\overline{1}n)-(1n)+(\overline{2}n)-(2n) +2\cdot(1n)+2\cdot(2n)}{2} =\\
&= \frac{Ch(1)+Ch(2)}{2} = \frac{(12)-(1\overline{2})+(12)+(1\overline{2})}{2} = (12).
\end{align*}
(The  division by $2$ is  eligible for degree one elements.)

In the general case, the proof goes by induction on $|H|$. For  $H~=~\{1,2,\dots, h\}$, denote $$\Sigma(h):=\sum_{S \subsetneq H} (Sn) * (-Ch(n))^{*\big(|H-S| -1\big)}.$$
In this sum, we first group together the summands such with $h \notin S$, and after those with $h \in S$.
Next, we apply the inductive assumption.
$$\Sigma(h)=
-Ch(n) * \Sigma(h-1)+(12\dots h-1\;n)+ (hn)*\Sigma(h-1) =$$
$$ ((n\overline{h})-(nh))*(12...h-1) +(nh)*(1...h-1)=
(12\dots h).$$

\medskip
Let us exemplify the induction step for $|H|=3$:
\begin{align*}
-&Ch(n) * \big(\underbrace{-Ch(n) + (1n)+(2n)}_{=(12)}\big)+(12n)+\underbrace{(23n)+(13n)+(3n)*(-Ch(n))}_{= \big((2n)+(1n)-Ch(n)\big) * (3n)} =\\
=& (12)*\big( \underbrace{-Ch(n) + (1n) + (3n)}_{=(23)} \big) =\\
=& (12)*(13) = (123).\qed
\end{align*}

\medskip

\section{The ring of  perfect manifolds. The second presentation}\label{SecPerf}

A nice manifold is \textit{perfect} if its label has no overlines.  That is, a perfect manifold is characterized by codirected collections of its edges.
\begin{prop}\label{PropPerf}
 As a $\mathbb{Z}$-module, the ring  $H^*(\mathcal{M}_3(L),\mathbb{Z})$  is generated by perfect manifolds.
\end{prop}
\textit{Proof.} The group $H^2(\mathcal{M}_3(L),\mathbb{Z})$ is generated by $(ij)$ and $(i\overline{j})$. Due to Propositions \ref{Ch} and \ref{Ch2}, each $(i\overline{j})$ is a linear combination of perfect manifolds. Indeed,
${(i\overline{j}) = (jk) - (ik) \text{ for any }k \neq i,j.}$

Now  let us show that any product of perfect manifolds  is a perfect manifold.

Let us call the perfect manifolds of type $(ij)$ \textit{elementary perfect manifolds}. Clearly, each perfect manifold is a product of elementary ones.

Now, we prove that the product of two perfect manifolds is an integer linear combination of   perfect manifolds.
We may assume that the second factor is an elementary perfect manifold, say, $(12)$.
Let the first factor be $(I_1)\cdot(I_2)\cdot(I_3)\cdot\ldots\cdot(I_k)$.

We need the following case analysis:\begin{enumerate}
                                      \item If at least one of $1,2$ does not belong to $\bigcup I_i$, the product is a perfect manifold by the  computation  rule (4).
                                      \item If $1$ and $2$ belong to different $I_i$, we use the {following:}

{for any perfect manifold $(I_1)\cdot(I_2)$ with $ i\in I_1, j\in I_2$, we have} $$(I_1)\cdot(I_2)\smile (ij)=(I_1)\smile (ij)\smile (I_2)= (I_1j)\smile (I_2)= (I_1\cup I_2).$$

                                      \item Finally, assume that $1,2 \in I_1$. Choose $i\notin I_1,\ \ j\notin I_1$ such that $i$ and $j$ do not belong to one and the same $I_k$.
By Corollary \ref{CorChern2}, $$(I_1)\cdot(I_2)\cdot(I_3)\cdot\dotso\cdot(I_k)\smile (12)=(I_1)\cdot(I_2)\cdot(I_3)\cdot\ldots\cdot(I_k)\smile \big((1i)+(2j)-(ij)\big).$$
After {expanding} the brackets, one reduces this to the above cases.
\qed
                                    \end{enumerate}

\bigskip

We are ready to describe a  presentation of the cohomology ring in terms of perfect manifolds.
Define the commutative graded ring $\mathcal{P}$  with the ring operation "$*$" as follows:

(\textbf{Generators}) {The ring $\mathcal{P}$ is a free $\mathbb{Z}$--module generated by }the labels of all existing perfect manifolds.
That is, the generators bijectively correspond to  all unordered formal products $$(I_1)\cdot ...\cdot(I_k),$$ where $I_1,...,I_k$ are some non-empty pairwise disjoint subsets of $[n]$.

(\textbf{Relations}) The relations in $\mathcal{P}$ are\begin{enumerate}
                                      \item  Let $I$ and $J$ be disjoint subsets of $[n]$. Then
                                      \begin{enumerate}

  \item $(I)* (J)=(I)\cdot (J).$

  \item $(Ii)* (Ji)=(IJi).$

\end{enumerate}
                                      \item  Length constraints:
                                        $(I)=0$  whenever the set $I$  is long.

                                      \item  The four-term relations:

$(ij)+(kl)=(jk)+(il)$ holds for any distinct $i,j,k,l$.

                                     This is equivalent to: \newline The value of  $(ij) +(ik)-(jk)$  does not depend on $j$ and $k$.

                                    \end{enumerate}

\bigskip

\textbf{Remark.}
All these relations hold true in $\widetilde{\mathcal{N}}$.

\begin{thm}
  The rings $\mathcal{P}$ and $H^*(\mathcal{M}_3(L),\mathbb{Z})$ are isomorphic.
\end{thm}
\textit{Proof.} We construct two ring homomorphisms $$\phi: \mathcal{P}\rightarrow \widetilde{\mathcal{N}}$$
and
$$\psi: \widetilde{\mathcal{N}}\rightarrow \mathcal{P}$$

and prove that (a)  $\phi$ is epimorphism, and  (b) $\psi \circ \phi=id$.

By definition, $\phi$ sends each $(I_1)\cdot ... \cdot (I_k)\in \mathcal{P}$ to an element of $\widetilde{\mathcal{N}}$ with the same label. The homomorphism $\phi$ is an
epimorphism. The proof goes analogously to that of Proposition  \ref{PropPerf}.

Now define $\psi$:
\begin{enumerate}
  \item Set   $\psi(i\overline{j}):=(jk)-(ik)$  for any $k\neq i$  and $k\neq j$. The target expression does not depend on $k$  due to the four-term relation.
  \item Set $\psi(I\overline{J}):=(I)\cdot(Jk)-(Ik)\cdot (J)$  for any $k \notin I$  and $k \notin J$. The target expression does not depend on $k$.
Indeed, $(I)\cdot (Jk)-(Ik)\cdot(J)= (I)\cdot(J)*\Big((jk)-(ik)\Big)$  for $i\in I, j \in J$. The four-term  relation $(jk)-(ik)=(jk')-(ik')$ completes the proof of correctness.
  \item Set  $\psi\Big((I_1\overline{J}_1)\cdot ...\cdot (I_m\overline{J}_m)\Big)$  as the product $\psi(I_1\overline{J}_1)* ...* \psi(I_m\overline{J}_m)$ .
\end{enumerate}

To prove the correctness,
we need to show that $\psi$ sends all the relations in $\widetilde{\mathcal{N}}$ to zero.

\begin{enumerate}
\item $\psi\Big((I\overline{J})+(J\overline{I})\Big) =0.$
This follows directly from the definition of the map $\psi$.

  \item For $I_1\cap I_2=\{i\},\ \ I_1\cap J_2= \emptyset, \ \ I_2\cap J_1= \emptyset,\  \hbox{and} \ \ J_1\cap J_2= \emptyset$, we need to prove that
$$\psi(I_1\overline{J}_1)* \psi(I_2\overline{J}_2)=\psi(I_1\cup I_2\overline{J_1\cup J_2}).$$

First, we prove that $\psi(i\overline{j})* \psi(i\overline{k})=\psi(i\overline{jk})$ by a straightforward calculation:
$$\psi(i\overline{j})* \psi(i\overline{k}) = \Big((jl)-(il)\Big)*\Big((jk)-(ij)\Big) = (kj)*\Big((jl)-(il)\Big) = \psi(i\overline{jk}).$$

Then the general case follows from the calculation:
\begin{align*}
\psi(I_1\overline{J}_1)* \psi(I_2\overline{J}_2)=& (-1)^{2} \cdot \psi(I_1)*\psi(i\overline{j})*\psi(\overline{J}_1)* \psi(I_2)*\psi(i\overline{k})*\psi(\overline{J}_2) =\\
=& \psi(I_1 \cup I_2)*\psi(\overline{J_1})*\psi(\overline{J_2})*\psi(i\overline{jk}) =\\
=& (-1)^{3}\cdot\psi(I_1 \cup I_2)*\psi(J_1)*\psi(J_2)*\psi(jk\overline{i}) =\\
=& (-1) \cdot\psi(I_1 \cup I_2)*\psi(J_{1}\cup J_{2} \overline{i}) =\\
=& (-1)^{2}\cdot\psi(I_1 \cup I_2)*\psi(i\; \overline{J_{1}\cup J_{2}}) =\psi(I_1\cup I_2\overline{J_1\cup J_2}).
\end{align*}

                                               \item  $\psi(I\overline{J})=0$  whenever one of the sets $I$ or $J$ is long. This follows from the definition.
                                               \item  $\psi\Big((ij)*(i\overline{j})\Big) =0.$ Indeed,  $\psi\Big((ij)*(i\overline{j})\Big) = (ij)*\Big((ju)-(iu)\Big) = (iju)- (iju)=0$.

                                      \item For any $i,j,k \in [n]$ such that $i\neq j, i\neq k$, $$\psi \Big((ij)-(i\overline{j})-(ik)+(i\overline{k})\Big)=0.$$

Indeed, $\psi \Big((ij)-(i\overline{j})-(ik)+(i\overline{k})\Big)= (ij) - (jk) + (ik)- (ik) +(kj)- (ij) =  0$.

                                    \end{enumerate}

\end{document}